\documentclass[11pt]{amsart}
\usepackage{amssymb}
\usepackage[mathscr]{eucal}

\theoremstyle{plain}
\newtheorem*{theorem}{Theorem}
\newtheorem*{proposition}{Proposition}
\newtheorem*{lemma}{Lemma}
\newtheorem*{corollary}{Corollary}

\theoremstyle{remark}

\def\C{\mathbb C}
\def\FH{F(H)}
\def\l{\lambda}
\def\la{\langle}
\def\ot{\otimes}
\def\PH{P(H)}
\def\SPH{SP(H)}
\def\R{\mathbb R}
\def\ra{\rangle}
\def\rank{\operatorname{rank}}
\def\tr{\operatorname{tr}}

\begin{document}
\title[]{Local automorphisms of some quantum mechanical structures}
\author{LAJOS MOLN\'AR}
\address{Institute of Mathematics and Informatics\\
         Lajos Kossuth University\\
         4010 Debrecen, P.O.Box 12, Hungary}
\email{molnarl@math.klte.hu}
\thanks{  This research was supported from the following sources:\\
          1) Joint Hungarian-Slovene research project supported
          by OMFB in Hungary and the Ministry of Science and
          Technology in Slovenia, Reg. No. SLO-2/96,\\
          2) Hungarian National Foundation for Scientific Research
          (OTKA), Grant No. T--030082 F--019322,\\
          3) A grant from the Ministry of Education, Hungary, Reg.
          No. FKFP 0304/1997}
%\subjclass{Primary: 47B49, 03G12}
%\keywords{Posets of idempotents, Jordan ring of
%selfadjoint operators, automorphisms, local automorphisms}
\date{\today}
\begin{abstract}
Let $H$ be a separable infinite dimensional complex Hilbert space.
We prove that every continuous 2-local automorphism of the poset
(that is, partially ordered set) of all idempotents on $H$ is an
automorphism.
Similar results concerning the orthomodular poset of all projections
and the Jordan ring of all selfadjoint operators on $H$ without the
assumption on continuity are also presented.
\end{abstract}
\maketitle

\newpage

\vskip 1cm

The orthomodular lattice or quantum logic of projections on a Hilbert
space plays
fundamental role in the mathematical foundations of quantum mechanics.
The interest in the poset of all skew projections (that is, idempotents)
has been also aroused (see \cite{Ovc} and the references therein) since
it can be defined on an arbitrary Banach space (or, more generally,
topological vector space).

In relation to any algebraic structure, the importance of the study of
automorphisms needs no justification. In a series of papers (see
\cite{MolStud1, MolGyor} and their references),
motivated by a problem of Larson \cite[Some concluding remarks (5),
p. 298]{Lar} (also see \cite{LarSour}) we investigated
the surprising and sometimes probably unexpected phenomenon
when the local automorphisms of a given Banach algebra are
all automorphisms. If this is the case, then one can say that
the local actions of the automorphism group
determines that group completely. In our papers we showed several Banach
algebras (in fact, they were mainly $C^*$-algebras) which possess this
property.
The aim of this present paper is to study a similar problem for posets
of idempotents and for the Jordan ring of selfadjoint operators.

Let us fix the notation and terminology that we shall use throughout.
Let $H$ be a Hilbert space. Let $B(H)$
stand for the set of all bounded linear operators on $H$. The set of all
skew projections (that is, idempotents) in $B(H)$ is denoted by $\SPH$,
and $\PH$ stands for the set of all projections.
For any $P,Q \in \SPH$ we write $P\leq Q$ if $PQ=QP=P$ and
we say that $P,Q$ are orthogonal if $PQ=QP=0$.
In what follows $\SPH$ is regarded as a poset with the relation $\leq$
and $\PH$ is viewed as an orthomodular poset with the additional map $P
\mapsto I-P$ of orthocomplementation \cite[Chapter 2]{Dvur}.

The transformation $\phi :\mathcal A \to \mathcal A$ of the algebraic
structure $\mathcal A$ is called a local automorphism
if for every $x\in \mathcal A$, there is an
automorphism $\phi_{x}$ of $\mathcal A$ for which
$\phi(x)=\phi_{x}(x)$.
The map $\phi$ is called a 2-local automorphism
if for every $x,y\in \mathcal A$, there is an
automorphism $\phi_{x,y}$ of $\mathcal A$ for which
$\phi(x)=\phi_{x,y}(x)$ and
$\phi(y)=\phi_{x,y}(y)$.
In our mentioned papers we considered linear maps on various
algebras which were local automorphisms and showed that they are in fact
automorphisms. The main result of \cite{BrSe} says that every linear map
which is a local automorphism of the algebra $B(H)$ on a separable
infinite dimensional Hilbert space $H$ is an automorphism (see
\cite{MolStud1} for a stronger result).
It turned out in \cite{Sem} that the assumption of linearity in
the result of \cite{BrSe}
can be dropped provided we pay the price that we
consider 2-local automorphisms instead of (1-)local automorphisms. More
precisely, the main result of \cite{Sem} tells us
that every 2-local automorphism of $B(H)$ (no linearity is assumed)
is an automorphism. For some results concerning 2-local automorphisms
of CSL algebras see \cite{Cri}.

It was proved in \cite{Ovc} that in the case of an infinite dimensional
separable Hilbert space $H$, the automorphisms of $\SPH$ with respect
to the partial order $\leq$ are exactly the maps
\[
P \longmapsto TPT^{-1}, \qquad
P \longmapsto TP^*T^{-1}
\]
where $T$ is an invertible bounded linear or conjugate-linear operator
of $H$.

\vskip .3cm
The main result of the paper is the following theorem.

\begin{theorem}
Let $H$ be an infinite dimensional separable complex Hilbert space.
Every continuous 2-local automorphism of the poset $\SPH$ is an automorphism.
\end{theorem}

It is easy to see that the corresponding assertion for
(1-)local automorphisms
fails to be true. In fact, the local automorphisms of $\SPH$ can
be almost arbitrary.

In the proof of our result we shall use the following easy lemma.

\begin{lemma}
Let $A,B \in B(H)$. Suppose that for every $x\in H$ we have either $\la
Ax,x\ra=0$ or $\la Bx,x \ra=0$. Then either $A=0$ or $B=0$.
\end{lemma}

\begin{proof}
The sets
$\{ x\in H \, : \, \la Ax,x\ra =0\}$ and
$\{ x\in H \, : \, \la Bx,x\ra =0\}$ are closed and their union is $H$.
By Baire's category theorem we obtain that one of these sets has
nonempty interior. So, we can suppose that there exist $x_0 \in H$ and
$\epsilon_0>0$ such that
\[
\la A(x_0+\epsilon x),x_0+\epsilon x \ra =0
\]
for every unit vector $x\in H$ and $0\leq \epsilon
<\epsilon_0$. This gives us that
\[
\la Ax_0,\epsilon x \ra + \la A(\epsilon x), x_0\ra +
\la A(\epsilon x),\epsilon x \ra =0,
\]
that is,
\[
\epsilon \la Ax_0, x \ra + \epsilon \la A x, x_0\ra +
\epsilon^2 \la Ax, x \ra =0.
\]
By the arbitrariness of $\epsilon$ we obtain that $\la Ax,x\ra =0$ for
every unit vector $x\in H$. This results in $A=0$.
\end{proof}

\begin{proof}[Proof of Theorem]
Let $\phi :\SPH\to \SPH$ be a continuous 2-local automorphism. Clearly,
$\phi$
preserves the partial order $\leq$ and the orthogonality between the
elements of $\SPH$. Furthermore, $\phi$ preserves the rank of the finite
rank idempotents.

In what follows we describe the form of $\phi$ on the set of all finite
rank elements of $\SPH$.
We first show that $\phi$ is finitely orthoadditive on this set in
question.
To see this, let $P,Q\in \SPH$ be orthogonal
and of finite rank. By the monotonity of $\phi$ we have
$\phi(P), \phi(Q) \leq \phi(P+Q)$. Since
$\phi(P)\phi(Q)=\phi(Q)\phi(P)=0$, it
follows that $\phi(P)+ \phi(Q) \leq \phi(P+Q)$.
By the rank preserving property of $\phi$ we obtain that
\[
\rank(\phi(P)+\phi(Q))=
\rank(\phi(P))+\rank(\phi(Q))=
\]
\[
\rank(P)+\rank(Q)=
\rank(P+Q)=
\rank(\phi(P+Q))
\]
which yields that
\begin{equation}\label{E:2hat}
\phi(P)+ \phi(Q) =\phi(P+Q).
\end{equation}

We now extend $\phi$ from the set of all finite rank
projections to a Jordan
homomorphism of the ideal $\FH$ of all finite rank operators in $B(H)$.
The additive map $J:\mathcal A \to \mathcal B$ between the rings
$\mathcal A, \mathcal B$ is called
a Jordan homomorphism if
\begin{equation}\label{E:vanjor}
J(xy+yx)=J(x)J(y)+J(y)J(x) \qquad (x,y\in \mathcal A).
\end{equation}
In $\mathcal A,\mathcal B$ are algebras, then besides additivity, the
Jordan homomorphisms are also supposed to be linear.
Clearly, every homomorphism and every antihomomorphism
are Jordan homomorphisms.

Let $H_d$ denote an arbitrary $d$-dimensional ($d$ is finite) subspace
of $H$.
Consider the natural embedding $B(H_d) \hookrightarrow B(H)$ and for any
$h\in H$ let $\phi_h$ be defined by $\phi_{h}(P)=\la \phi(P)h, h
\ra$ for every projection $P$ on $H_d$.
Since $\phi$ is continuous, thus $\phi_h$ is a bounded orthoadditive
function. If $d\geq 3$, then by Gleason's theorem \cite[Theorem
3.2.16]{Dvur} there exists an operator $T_{h}$ on $H_d$ such that
\begin{equation}\label{E:2gleas}
\phi_{h} (P)=\tr T_{h} P \qquad (P \in P(H_d)).
\end{equation}
Let $P_1, \ldots, P_n\in \PH$ be finite rank projections (their
pairwise orthogonality is not assumed) and let $\l_1, \ldots , \l_n$
be complex numbers. Define
\begin{equation}\label{E:2het}
\psi(\sum_k \l_k P_k)=\sum_k \l_k \phi(P_k).
\end{equation}
We have to check that $\psi$ is well-defined. To see this, let
$P'_1, \ldots, P'_n\in \PH$ be finite rank projections and
$\mu_1, \ldots , \mu_n \in \C$ be such that
\[
\sum_k \l_k P_k=\sum_k \mu_k P'_k.
\]
Let $H_d$ be a finite dimensional subspace of $H$ of dimension $d\geq
3$ such that for the orthogonal projection $P_{H_d}$ onto $H_d$ we have
$P_k, P'_k \leq P_{H_d}$ $(k=1,\ldots ,n)$.
Let $T_{h}$ denote the linear operator on $H_d$
corresponding to $\phi_{h}$ (see \eqref{E:2gleas}). We compute
\[
\la \sum_k \l_k \phi(P_k)h,h\ra =
\sum_k \l_k \phi_h(P_k) =
\tr T_{h} (\sum_k \l_k P_k)=
\]
\[
\tr T_h (\sum_k \mu_k P'_k)=
\sum_k \mu_k \phi_h(P'_k) =
\la \sum_k \mu_k \phi(P'_k)h,h\ra .
\]
Since this holds true for every $h\in H$, we obtain that $\psi$ is
well-defined. As $F(H)$ is the linear span of its projections, the
definition \eqref{E:2het}
clearly implies that $\psi$ is a linear transformation on $F(H)$.
Since $\psi$ sends projections to idempotents,
it is now a standard argument to verify that $\psi$ is a Jordan
homomorphism of $F(H)$.
See, for example, the proof of \cite[Theorem 2]{MolStud1}.

As $F(H)$ is a locally matrix ring, it follows from a classical
result
of Jacobson and Rickart \cite[Theorem 8]{JR} that $ \psi$ can
be written as $ \psi= \psi_1 +\psi_2$, where $\psi_1$ is a
homomorphism and $\psi_2$ is an antihomomorphism. Let $P\in \PH$ be
rank-one. Because $ \psi(P)=\phi(P)$ is also rank-one, we obtain that
one of the idempotents $\psi_1(P), \psi_2(P)$ is zero. Since $F(H)$ is a
simple ring, it is easy to see that this implies that either
$\psi_1$ or $\psi_2$ is identically zero, that is, $ \psi$ is
either a homomorphism or an antihomomorphism of $F(H)$. In what follows
we can assume without loss of generality that $ \psi$ is a
homomorphism.

We show that $ \psi$ preserves the rank. Let $A\in F(H)$ be a
rank-$n$ operator. Then there is a rank-$n$ projection $P$ such that
$PA=A$. The rank of $\phi(P)$ is also $n$.
We have $ \psi(A)= \psi(P) \psi(A)=\phi(P) \psi(A)$
which proves that $ \psi(A)$ is of rank at most $n$. If $Q$ is
any rank-$n$ projection, then there are finite rank operators $U,V$
such that $Q=UAV$. Since $\phi(Q)= \psi(Q)= \psi(U)
\psi (A)
\psi(V)$ and the rank of $\phi(Q)$ is $n$, it follows that the rank of
$ \psi(A)$ is at least $n$. Therefore, $ \psi$ is rank
preserving. We now refer to Hou's work \cite{Hou} on the form of linear
rank preservers on operator algebras.
It follows from the argument leading to \cite[Theorem 1.2]{Hou}
that there are linear
operators $T,S$ on $H$ such that $ \psi$ is of the form
\begin{equation}\label{E:2gleasegy}
\psi(x\ot y)=(Tx)\ot (Sy) \qquad (x,y \in H)
\end{equation}
(recall that we have assumed that $ \psi$ is a homomorphism).
Here, for any $u,v \in H$, $u\ot v$ stands for the operator defined
by $(u\ot v)(z)=\la z,v\ra u$ $(z\in H)$.
We claim that $T,S$ are bounded. This follows from \cite[Lemma
1]{MolLAA} which states that
if $T,S$ are linear operators on $H$ with the property that the map $x
\mapsto (Tx) \ot (Sx)$ is continuous on the unit ball of $H$, then
$T, S$ are bounded.
We infer from
\eqref{E:2gleasegy} that $\la Tx, Sx\ra=\la x,x\ra$ for every
unit vector $x\in H$ ($\phi$ sends rank-one
projections to idempotents).
By polarization this implies that $\la Tx, Sy\ra=\la x,y\ra$ $(x,y\in
H)$.
Consequently, we have $S^*T=I$. From \eqref{E:2gleasegy} we deduce that
$\phi(P)=TPS^*$ for every rank-one
projection $P$. By the finite orthoadditivity of $\phi$
appearing in \eqref{E:2hat}, it follows that $\phi(P)=TPS^*$
holds true for every finite-rank projection $P$ as well.

Since $S^*T=I$, it follows that $Q=TS^*$ is an idempotent. We can write
\[
\phi(P)=Q\phi(P)Q+Q\phi(P)(I-Q)+(I-Q)\phi(P)Q+(I-Q)\phi(P)(I-Q)
\]
for every $P\in \SPH$. We claim that the two middle terms on the right
hand side of the equality above are in fact missing. Denote
\[
\phi_{11}(P)=Q\phi(P)Q, \quad
\phi_{12}(P)=Q\phi(P)(I-Q),
\]
\[
\phi_{21}(P)=(I-Q)\phi(P)Q, \quad
\phi_{22}(P)=(I-Q)\phi(P)(I-Q).
\]
Let $P\in \PH$ be
fixed and let $P'$ be an arbitrary finite rank projection with
$P'\leq P$. We know that $\phi(P') \leq \phi(P)$.
Since $\phi(P')=TP'S^*$, we obtain that $\phi(P')Q=Q\phi(P')=\phi(P')$.
Hence,
\[
\phi(P')\phi_{11}(P)=(\phi(P')Q)\phi(P)Q=
\phi(P')\phi(P)Q=
%\]
%\[
\phi(P')Q=
\phi(P').
\]
Therefore, $TP'S^*\phi_{11}(P)=TP'S^*$. Similarly, we have
$\phi_{11}(P)TP'S^*=TP'S^*$. By the arbitrariness of $P'$ it
follows that
\[
\phi_{11}(P)TPS^*=TPS^*\phi_{11}(P)=TPS^*.
\]
This means that
\begin{equation}\label{E:2val}
TPS^*\leq \phi_{11}(P).
\end{equation}
The local property of $\phi$ implies that $\phi(P)+\phi(I-P)=I$.
Therefore, $\phi(I-P)=I-\phi(P)$.
Writing $I-P$ for $P$ in \eqref{E:2val}, we have
\[
Q-TPS^*=T(I-P)S^*\leq \phi_{11}(I-P)=Q-\phi_{11}(P).
\]
We deduce that $TPS^*=\phi_{11}(P)$.
We next compute
\[
TP'S^*\phi_{12}(P)=
\phi(P')Q\phi(P)(I-Q)=
\phi(P')\phi(P)(I-Q)=
\]
\[
\phi(P')(I-Q)=
\phi(P')Q(I-Q)=0.
\]
Since this holds for every finite rank projection $P'$ for which $P'\leq
P$, we infer that $TPS^*\phi_{12}(P)=0$. Since $\phi(I-P)=I-\phi(P)$, we
have $\phi_{12}(I-P)=-\phi_{12}(P)$. Therefore,
\[
0=T(I-P)S^*\phi_{12}(I-P)=
TS^*(-\phi_{12}(P))-TPS^*(-\phi_{12}(P))=
\]
\[
-Q\phi_{12}(P)+TPS^*\phi_{12}(P)=
-\phi_{12}(P)+TPS^*\phi_{12}(P)=
-\phi_{12}(P)
\]
and, hence,
we obtain $\phi_{12}=0$. Similarly, one can verify that $\phi_{21}=0$.
Consequently, our map $\phi$ is of the form
\begin{equation}\label{E:2to}
\phi(P)=TPS^*+\phi_{22}(P) \qquad (P\in \PH ).
\end{equation}
We know that $\phi_{22}(P)=0$ for every finite rank
projection $P$ (recall that $\phi(P)=TPS^*$ for all such $P$). We claim
that $\phi_{22}(P)=0$ holds for every $P\in \SPH$ as well.

Assume for a moment that
$\phi_{22}(P)\neq 0$ for every projection $P$ of infinite rank and
infinite corank.
We can choose uncountably many projections $P_\alpha$ of infinite rank
and infinite corank such that $P_\alpha P_\beta $ is a finite rank
projection for every $\alpha \neq \beta$ (see, for example, the
proof of \cite[Theorem 1]{MolStud1}).
Using the local form of $\phi$ we see that the rank of
the idempotent $\phi(P_\alpha)\phi(P_\beta)$ is equal to the rank of
$P_\alpha P_\beta$. On the other hand, referring to
the injectivity of $T$ and to the surjectivity of $S^*$ (these follow
from $S^*T=I$), we find that the
rank of $TP_\alpha S^* TP_\beta S^*=TP_\alpha P_\beta S^*$ is the same
as that of $P_\alpha P_\beta$.
By \eqref{E:2to} this gives us that
the rank of
$\phi_{22}(P_\alpha)\phi_{22}(P_\beta)$ is 0, that is, we have
\[
\phi_{22}(P_\alpha)\phi_{22}(P_\beta)=0
\]
for every $\alpha\neq \beta$. This means that the range of $\phi_{22}$
contains uncountably many nonzero pairwise orthogonal idempotents which
plainly contradicts the separability of $H$.
Therefore,
$\phi_{22}(P)=0$ holds for a projection $P$ of infinite rank and
infinite corank. The projections $P$ and
$I-P$ can be connected by a continuous curve inside the set of
projections (this is an easy consequence of the arcwise connectedness
of the unitary group of $B(H)$). If $R, R'\in \PH$ are lying on the
same arc in $\PH$, then by the continuity of
$\phi$, the idempotents $\phi_{22}(R)$, $\phi_{22}(R')$ are close enough
to each other if $\| R-R'\|$ is sufficiently small.
Taking into account that the norm of a nonzero
idempotent is not less than 1, this, together with $\phi_{22}(P)=0$,
yields that $\phi_{22}(I-P)=0$. Hence, we have
$I-Q=\phi_{22}(I)=\phi_{22}(P)+\phi_{22}(I-P)=0$
and thus $\phi_{22}=0$.
Since $TS^*=Q=I=S^*T$, we obtain $S^*=T^{-1}$. Therefore, $\phi$ is of
the form
\[
\phi(P)=TPT^{-1} \qquad (P\in \PH)
\]
where $T$ is an invertible bounded linear operator
on $H$. We show that our map $\phi$ is of this form on the whole
set $\SPH$.
To verify this, we can obviously suppose that $T=I$. So, assume
that $\phi$ is the identity on the set of all projections.

Let $P$ be any idempotent. Pick an arbitrary unit vector $x\in H$ and
consider the operator $\phi(x\ot x)\phi(P)\phi(x\ot x)$. Taking into
account the local property of $\phi$, the form of the automorphisms of
$\SPH$ and that $\phi$ is the identity on the set of all projections,
we have either
\[
\la \phi(P)x,x\ra x\ot x=
\phi(x\ot x)\phi(P)\phi(x\ot x)=
\]
\[
A \cdot x\ot x \cdot A^{-1} APA^{-1} A \cdot x\ot x \cdot A^{-1}=
A\cdot (\la Px,x \ra x\ot x) \cdot A^{-1}=
\]
\[
\la Px,x\ra \phi(x\ot x)=
\la Px,x\ra x\ot x,
\]
or, using a similar computation,
\[
\la \phi(P)x,x\ra x\ot x=
\overline{\la Px,x\ra} x\ot x,
\]
or
\[
\la \phi(P)x,x\ra x\ot x=
\la P^*x,x\ra x\ot x,
\]
or
\[
\la \phi(P)x,x\ra x\ot x=
\overline{\la P^*x,x\ra} x\ot x.
\]
This gives us that for every $x\in H$ we have either $\la \phi(P)x,x\ra
=\la
Px,x \ra$ or $\la \phi(P)x,x\ra =\la P^*x,x \ra$. Our lemma implies that
for any $P\in \SPH$ we have either $\phi(P)=P$ or $\phi(P)=P^*$.
We assert that this results in either
$\phi(P)=P$ for all $P\in \SPH$ or
$\phi(P)=P^*$ for all $P\in \SPH$. Indeed, let $P\in \SPH$ be a
non-selfadjoint
finite rank idempotent for which $\phi(P)=P$. Consider any
non-selfadjoint finite rank
idempotent $P'$ which is orthogonal to $P$. If $\phi(P')={P'}^*$, then
by $\phi(P)+\phi(P')=\phi(P+P')$ we would arrive at a contradiction.
So, $\phi(P')=P'$ for every finite rank idempotent which is orthogonal
to $P$. This implies that $\phi(R)=R$ for every finite rank
idempotent with $P\leq R$. Let $P'$ be any non-selfadjoint finite rank
idempotent. Suppose that $\phi(P')={P'}^*$. If $R$ is any
non-selfadjoint finite rank idempotent for
which $P'\leq R$, then we have similarly as before that $\phi(R)=R^*$.
Now, if $P,P' \leq R$, then we obtain on the one hand that $\phi(R)=R$
and
on the other hand that $\phi(R)=R^*$. But this is a
contradiction. Therefore,
we have $\phi(P')=P'$ for every finite rank idempotent $P'$.
By the monotonity of $\phi$ it now follows that $P\leq \phi(P)$
for every $P\in \SPH$.
Putting $I-P$ in the place of $P$ we finally obtain $\phi(P)=P$ $(P\in
\SPH)$.

In case $ \psi$ is an antihomomorphism, we can follow a similar
argument. The proof is complete.
\end{proof}

Examining the proof of our theorem, we can reach a result of the
same spirit concerning the orthomodular poset $\PH$. First we need
the form of the automorphisms of $\PH$ as an orthomodular poset
(as it was mentioned in \cite[Remark 4.5.]{Ovc}, an automorphism
of $\PH$ as a poset is not necessarily an automorphism of it as an
orthomodular poset). So, let $\phi :\PH \to \PH$ be a bijection
which preserves the partial order $\leq$ as well as the
orthocomplementation $P\mapsto I-P$ in both directions. We easily
get that $\phi$ preserves the orthogonality between the elements
of $\PH$ and then that $\phi$ is orthoadditive. Using, for
example, the result of Bunce and Wright \cite{BW} which solves the
Mackey-Gleason problem (also see \cite{BW2, Dye}), we see that
$\phi$ can be extended to a bounded linear transformation of
$B(H)$. Since this map sends projections to projections, one can
verify that it is a Jordan *-homomorphism of $B(H)$. Since its
range contains a rank-one operator and an operator with dense
range (in fact, the range contains every projection), it is a
trivial consequence of \cite[Theorem 1]{MolStud1} that this
transformation is a Jordan *-automorphism of $B(H)$ onto itself.
As $B(H)$ is a prime algebra (that is, $AB(H)B=\{ 0\}$ implies
that either $A=0$ or $B=0$), it follows from a classical theorem
of Herstein \cite[Theorem H]{Her} that our map is either a
*-automorphism or a *-antiautomorphism. The structure of those
maps is well-known. In fact, they are of the form
\[
A\longmapsto UAU^* , \qquad
A\longmapsto VA^*V^*
\]
where $U$ is unitary and $V$ is antiunitary on $H$. So, $\phi$ is either
of the form
\[
\phi(P)=UPU^* \qquad (P\in \PH)
\]
or of the form
\[
\phi(P)=VPV^* \qquad (P\in \PH)
\]
where $U$ is
unitary and $V$ is antiunitary on $H$. Now,
we can prove the following result.

\begin{proposition}
Every 2-local automorphism of the orthomodular poset $\PH$ is an
automorphism.
\end{proposition}

\begin{proof}
One can just follow the proof of our theorem.
The only thing which deserves checking is that here we do
not need
the continuity of $\phi$. To see this, let us go through those parts of
the previous proof where we have used continuity.
The first such place was where we applied Gleason's theorem. But as
$\phi$ sends projections to projections, it follows that $\phi$ is
bounded, so we do not need continuity here.
We next used the continuity when showing that the operators $T,S$ in
\eqref{E:2gleasegy} are continuous. Since in our present case $\psi$
sends
projections to projections, it follows that for every unit vector $x\in
H$, the operator $Tx \ot Sx$ is a projection. This implies that $Tx=Sx$
and $\|Tx\|=\| Sx\|=1$. So, $T=S$ is an isometry.
The third appearence of the continuity of $\phi$ was where we proved
that if $\phi_{22}(P)=0$ for
a projection $P$ of infinite rank and infinite corank, then
$\phi_{22}=0$. In our present case the terms in the
decomposition $\phi(R)=\phi_{11}(R)+\phi_{22}(R)$ $(R\in \PH)$ are
projections. Let $P'\in \PH$ be such that $\| P-P'\|<1$. By the
2-local property of $\phi$ we see that $\|\phi(P)-\phi(P')\| =\| P-P'\|
<1$. Since
\[
\| \phi_{22}(P')\|=
\| \phi_{22} (P)-\phi_{22}(P')\| \leq
\|\phi(P)-\phi(P')\| <1,
\]
we deduce that $\phi_{22}(P')=0$. As we can go from $P$ to $I-P$ in
finitely many steps $P=P_0, P_1, \ldots, P_n=I-P$ such that
$\| P_{k-1}- P_k\| <1$, we obtain that $\phi_{22}(I-P)=0$. This
gives us that $\phi_{22}(I)=0$ which implies $\phi_{22}=0$.
The proof is complete.
\end{proof}

Beside the various structures of projections and idempotents, the Jordan
ring $B(H)_h$ of all selfadjoint operators on $H$ is also well-known to
be of great importance in the mathematical description of quantum
mechanics. Our final result shows that the 2-local automorphisms of
$B(H)_h$ are necessarily automorphisms. To prove this, we first
describe the
automorphisms in question. So, let $\phi:B(H)_h \to B(H)_h$ be a Jordan
automorphism (this means that $\phi$ is an addivite bijection satisfying
\eqref{E:vanjor}).
We have
\[
\phi(ABA)=\phi(A)\phi(B)\phi(A)
\qquad (A,B \in B(H)_h)
\]
(see, for example, \cite[Lemma 2]{Her}). One can readily verify
that $\phi$ preserves the partial order and the orthogonality between
the
projections. Therefore, we obtain that $\phi$ sends rank-one projections
to rank-one projections, if $\{ P_n\}_n$ is a maximal system of pairwise
orthogonal rank-one projections then the same holds for $\{
\phi(P_n)\}_n$, and $\phi(I)=I$. Let $P$ be any rank-one projection.
Since $\phi(P)$ is also rank-one, if $\lambda \in \R$, then from
the equality
\[
\phi(\lambda P)=\phi(P(\lambda P)P)=\phi(P)\phi(\lambda P)\phi(P)
\]
it follows that $\phi(\lambda P)=f(\lambda )\phi(P)$ for some real
number $f(\lambda)$. It is easy to verify that $f :\R \to \R$ is a
ring-homomorphism with $f(1)=1$. It is well-known and, in fact, it
requires only elementary analysis
to prove that this implies that $f(\lambda)=\lambda$ $(\lambda \in \R)$.
So, we have $\phi(\lambda P)=\lambda \phi(P)$. Since
\[
2\lambda \phi(P)=\phi(2\lambda P)=\phi(\lambda
I)\phi(P)+\phi(P)\phi(\lambda I)
\]
for every rank-one projection, choosing a maximal orthogonal family of
such projections, we obtain that
\[
2\lambda I=\phi(\lambda I)I+I\phi(\lambda I).
\]
That is, we have
\[
\phi(\lambda I)=\lambda I \qquad (\lambda \in \R).
\]
It is now apparent that $\phi$ is real-linear. Let us define $\tilde
\phi :B(H) \to B(H)$ by
\[
\tilde \phi (A+iB)=\phi(A)+i\phi(B) \qquad (A,B \in B(H)_h).
\]
It is easy to check that $\tilde \phi$ is a (linear) Jordan
*-automorphism of $B(H)$. Hence, just as in the discussion right before
the formulation of our proposition, we find that $\phi$ is of the form
\[
\phi(A)=UAU^* \qquad (A\in B(H)_h)
\]
where $U$ is an either unitary or antiunitary operator on $H$.

Now, we are in a position to prove our final result which follows.

\begin{corollary}
Every 2-local automorphism of the Jordan ring $B(H)_h$ is an
automorphism.
\end{corollary}

\begin{proof}
Let $\phi: B(H)_h \to B(H)_h$ be a 2-local automorphism. Clearly,
$\phi_{|P(H)}$ is a 2-local automorphism of $P(H)$. By Proposition
we obtain that there exists a unitary or antiunitary operator $U$ on $H$
such that
\[
\phi(P)=UPU^* \qquad (P\in P(H)).
\]
We can assume without loss of generality that $\phi(P)=P$ for every
$P\in P(H)$. Now, similarly as in the proof of our theorem, picking any
$A\in B(H)_h$ and unit vector $x\in H$, considering the operator
$\phi(x\ot x)\phi(A)\phi(x\ot x)$, we find that
\[
\la \phi(A)x,x\ra x\ot x =\la Ax,x\ra x\ot x
\]
which implies that
\[
\la \phi(A)x,x\ra =\la Ax,x\ra .
\]
Therefore, we have $\phi(A)=A$ $(A\in B(H)_h)$ and this completes the
proof. \end{proof}

%Bibliography
\bibliographystyle{amsplain}

\begin{thebibliography}{99}

\bibitem{BrSe}
M. Bre\v sar and P. \v Semrl,
\emph{On local automorphisms and mappings that preserve idempotents},
Studia Math. \textbf{113} (1995), 101--108.

\bibitem{BW}
L.J. Bunce and D.M. Wright,
\emph{The Mackey-Gleason problem},
Bull. Amer. Math. Soc. \textbf{26} (1992), 288--293.

\bibitem{BW2}
L.J. Bunce and D.M. Wright, \emph{On Dye's theorem for Jordan
operator algebras}, Expo. Math. \textbf{11} (1993), 91--95.

\bibitem{Cri}
R. Crist, \emph{Local automorphisms}, Proc. Amer. Math. Soc.
\textbf{128} (2000), 1409--1415.

\bibitem{Dye}
H.A. Dye, \emph{On the geometry of projections in certain operator
algebras}, Ann. Math. \textbf{61} (1955), 73--89.

\bibitem{Dvur}
A. Dvure\v censkij,
\emph{Gleason's Theorem and Its Applications},
Kluwer Academic Publishers, 1993.

\bibitem{Her}
I.N. Herstein,
\emph{Jordan homomorphisms},
Trans. Amer. Math. Soc. \textbf{81} (1956), 331--341.

\bibitem{Hou}
J.C. Hou,
\emph{Rank-preserving linear maps on $B(X)$},
Sci. China Ser. A  \textbf{32} (1989), 929--940.

\bibitem{JR}
N. Jacobson and C. Rickart,
\emph{Jordan homomorphisms of rings},
Trans. Amer. Math. Soc. \textbf{69} (1950),  479--502.

\bibitem{Lar}
D.R. Larson,
\emph{Reflexivity, algebraic reflexivity and linear interpolation},
Amer. J. Math. \textbf{110} (1988), 283--299.

\bibitem{LarSour}
D.R. Larson and A.R. Sourour,
\emph{Local derivations and local automorphisms of $B(X)$},
in Proc. Sympos. Pure Math. 51, Part 2, Providence, Rhode
Island 1990, 187--194.

\bibitem{MolStud1}
L. Moln\'ar,
\emph{The set of automorphisms of $B(H)$ is topologically reflexive in
$B(B(H))$},
Studia Math. \textbf{122} (1997), 183--193.

\bibitem{MolLAA}
L. Moln\'ar,
\emph{Some multiplicative preservers on $B(H)$},
Linear Algebra Appl. \textbf{301} (1999), 1--13.

\bibitem{MolGyor}
L. Moln\'ar and M. Gy\H{o}ry,
\emph{Reflexivity of the automorphism and isometry groups of the
suspension of $B(H)$},
J. Funct. Anal. \textbf{159} (1998), 568--586.

\bibitem{Ovc}
P.G. Ovchinnikov,
\emph{Automorphisms of the poset of skew projections},
J. Funct. Anal. \textbf{115} (1993), 184--189.

\bibitem{Sem}
P. \v Semrl, \emph{Local automorphisms and derivations on $B(H)$},
Proc. Amer. Math. Soc. \textbf{125} (1997), 2677--2680.

\end{thebibliography}

\end{document}